\documentclass[12pt,reqno]{article}
\usepackage{amssymb, amsmath}
\usepackage{amssymb,amsthm}
\newtheorem{thm}{Theorem}

\numberwithin{Definition}{section}
\numberwithin{Lemma}{section}
\numberwithin{Proposition}{section}

\numberwithin{Corollary}{section}
\numberwithin{Example}{section}
\numberwithin{Remark}{section}
\begin{document}
\begin{center}
\textbf{ASYMPTOTIC EXPANSION OF THE WAVELET TRANSFORM
FOR SMALL $a$ }
\end{center}
\begin{center}
R S Pathak and Ashish Pathak
\end{center}
\begin{center}
\textbf{Abstract }
\end{center}
Asymptotic expansion of the wavelet transform for small values of the dilation parameter a is obtained using asymptotic expansion of the Mellin convolution technique of Wong. Asymptotic expansions of Morlet wavelet transform, Mexican hat wavelet transform and Haar wavelet transform are obtained as special cases. 
\footnote{ This work is contained in the research monograph " The Wavelet Transform" by Prof. R S Pathak and
edited by Prof. C. K. Chui (Stanford University, U.S.A.) and published by Atlantis Press/World Scientific
(2009), ISBN: 978-90-78677-26-0, pp:164-169. 
}
\section{Introduction}The wavelet transform of f with respect to the wavelet $ \psi $ is defined by
\begin{equation}
\label{eq:Int-310.2}
 (W_\psi f)(b,a) = \frac{1}{\sqrt{a}}\int_{-\infty}^\infty f(t)
  \overline{\psi\left(\frac{t-b}{a}\right)}\,dt,b\in \mathbb{R},a>0,
\end{equation}
provided the integral exists ~\cite{Debnath}.\\
 For fixed parameter $ b \in \mathbb{R}$ in (\ref{eq:Int-310.2}) the family $ \left\{\psi\left(\frac{\bullet-a}{b}\right):a>0\right\} $ zooms in an every detail of in a neighborhood of $b$ as long as a is sufficiently small $a$. The frequency resolution is controlled by the parameter $a$ and for small $a$, $ (W_\psi f)(b,a)$ represents the frequency components of the signal $f$. Therefore, it is highly desirable to know the asymptotic behavior of $ (W_\psi f)(b,a)$ for small values of $a$. Using Fourier transform (\ref{eq:Int-310.2}) can also be expressedas
 \begin{equation}
\label{eq:Int-311.2}
 (W_\psi f)(b,a) = \frac{\sqrt{a}}{2\pi}\int_{-\infty}^\infty e^{ib\omega}  \hat f(\omega)
  \overline{\hat \psi(a\omega)}d\omega,
\end{equation}
where
\begin{equation}
\nonumber \hat f(\omega)= \int_{-\infty}^\infty e^{-it\omega} f(t)dt.
\end{equation}
Putting $a=\frac{1}{c}$, from (\ref{eq:Int-311.2}) we have
\begin{equation}
\label{eq:Int-312.2}
 (W_\psi f)(b,a) = \frac{\sqrt{c}}{2\pi}\int_{-\infty}^\infty e^{ibcu}  \hat f(cu)
  \overline{\hat \psi(u)}du,
\end{equation}
Asymptotic expansion with explicit error term for the general integral
\begin{eqnarray}
\label{eq:Int-313.2}
 I(x) = \int_0^\infty g(t) h(xt) dt,
\end{eqnarray}
as $ x \rightarrow \infty+ $, was obtained by Wong ~\cite{Wong}, ~\cite{Rwong} under different conditions on $g$ and $h$. The asymptotic expansion for (\ref{eq:Int-312.2}) can be obtained by setting $ h(t)=e^{ibt}\hat{f}(t) $ for fixed $b \in \mathbb{R}$. Let us recall basic result from ~\cite{Rwong}
that will be used in the present investigation.\\
\,\,\, Here we assume that $g(t)$ has an expansion of the form
\begin{eqnarray}
\label{eq:Int-315.2}
 g(t)\sim \sum_{s=0}^{n-1} c_s  t^{s+\lambda-1} \,\,,\,\, as \,\,\,\, t\rightarrow 0+,
\end{eqnarray}
where $ 0<\lambda \leq 1 $. Regarding the function $h$, we assume that as $t\rightarrow 0+,$
\begin{eqnarray}
\label{eq:Int-316.2}
 h(t)=O(t^\rho), \,\,\, \rho+\lambda>0,
\end{eqnarray}
and that as $t\rightarrow + \infty $
\begin{eqnarray}
\label{eq:Int-317.2}
 h(t) \sim  exp\left(i \tau t^p \right) \sum_{s=0}^\infty b_s t^{-s-\beta},
\end{eqnarray}
where $\tau\neq 0 $ is real, $p\geq1$ and $\beta>0$. Let $M[h;z]$ denote the generalized Mellin transform of h defined by
\begin{eqnarray}
\label{eq:Int-318.2}
M[h;z]= \lim_{\epsilon \rightarrow 0+}\int^\infty_0 t^{z-1} h(t)\,\,  exp\left(-\epsilon \, t^p \right)dt.
\end{eqnarray}
This together with (\ref{eq:Int-313.2}) and ~\cite[p.217]{Rwong}, gives
\begin{eqnarray}
\label{eq:Int-319.2}
I(x) = \sum^{n-1}_{s=0}c_s M[h;
s+\lambda]x^{-s-\lambda}+\delta_n(x),
\end{eqnarray}
where
\begin{eqnarray}
\label{eq:Int-320.2}
\delta_n(x)=\lim_{\varepsilon \rightarrow
0+}\int^{\infty}_{0}g_n(t)h(xt)\exp(-\varepsilon t^p)dt.
\end{eqnarray}
If we now define recursively $h^\circ (t)=h(t)$  and
\begin{eqnarray*}
h^{(-j)}(t)=-\int^{\infty}_{t}h^{(-j+1)}(u) du, \;\;\; j=1,2,...,
\end{eqnarray*}
then conditions of validity of aforesaid results are given by the
following ~\cite[Theorem 6, p.217]{Rwong}.
\begin{thm}
\label{thm:ch1sec2-2.1}
 Assume that (i) $g^{(m)}(t)$ is continuous on $(0,\infty)$, where
$m$ is a non-negative integer; (ii) $g(t)$ has an
expansion of the form(\ref{eq:Int-315.2}), and the expansion is $m$
times differentiable; (iii) $h(t)$ satisfies
(\ref{eq:Int-316.2}) and (\ref{eq:Int-317.2}) and (iv) and as $t
\rightarrow \infty, t^{-\beta}g^{(j)}(t)=O(t^{-1-\varepsilon})$  for $j=0,
1,..., m$ and for some $\varepsilon > 0.$ Under
these conditions, the result(\ref{eq:Int-319.2}) holds with
\begin{eqnarray}
\label{eq:Int-321.2}
\delta_n(x)=\frac{(-1)^m}{x^m}\int^{\infty}_{0} g^{(m)}_{n}(t)
h^{(-m)}(xt) dt,
\end{eqnarray}
where $n$ is the smallest positive integer such
that $\lambda+n > m.$
\end{thm}
The asymptotic expansion for the wavelet transform (\ref{eq:Int-311.2}) for large values of dilation parameter a has already been obtained in ~\cite{Rpathak}.\\
The aim of the present paper is to derive asymptotic expansion of  the wavelet transform given by (\ref{eq:Int-311.2}) for small values of $a$, using formula (\ref{eq:Int-319.2}).
\section{ASYMPTOTIC EXPANSION FOR SMALL $a$}
In this section using aforesaid technique, we obtain asymptotic expansion of $\left(W_\psi f\right)(b,a)$ for small values of  $a$, keeping $b$ fixed. We have
\begin{eqnarray}
\label{eq:Int-322.2}
\left(W_\psi f\right)(b,a)
\nonumber & = & \frac{\sqrt{c}}{2\pi}\int^\infty_{-\infty} e^{ibcu}\overline{\hat\psi}(u)\hat f(cu) du \\
    \nonumber & = & \frac{\sqrt{c}}{2\pi}\biggr\{\int^\infty_0 e^{ibcu}\overline{\hat\psi}(u)\hat f(cu) du \\ \nonumber && + \;
    \int^\infty_0 e^{-ibcu}\overline{\hat\psi}(-u)\hat f(-cu) du \biggr\} \,\,\,\,\,\,\
    \\ &=&  \frac{\sqrt{c}}{2\pi} \left(I_1 + I_2 \right),
    \,\,\,\,\ (say).
\end{eqnarray}
Let us set
\begin{eqnarray}
\label{eq:Int-323.2}
h(u)= e^{ibu} \hat{f}(u).
\end{eqnarray}
Assume that
\begin{eqnarray*}
\hat{f}(u) \sim \sum_{r=0}^\infty b_r u^{-r-\beta}, \,\,\, \beta > 0, \,\,\, u\rightarrow \infty;
\end{eqnarray*}
so that
\begin{eqnarray}
\label{eq:Int-324.2}
h(u) \sim e^{ibu} \sum_{r=0}^\infty b_r u^{-r-\beta}, \,\,\, \beta > 0, \,\,\, u\rightarrow \infty, b\neq0.
\end{eqnarray}
\begin{eqnarray*}
\overline{\hat{\psi}}(u) \sim \sum_{s=0}^\infty c_s u^{s+\lambda-1}, \,\,\, as \,\,\, u\rightarrow 0.
\end{eqnarray*}
For $n\geq1$, we write
\begin{eqnarray}
\label{eq:Int-325.2}
\overline{\hat{\psi}}(u) \sim \sum_{s=0}^{n-1} c_s u^{s+\lambda-1}+\overline{\hat{\psi}_n}(u),
\end{eqnarray}
where $0<\lambda\leq 1$. Also assume that
\begin{eqnarray}
\label{eq:Int-326.2}
h(u)= O(u^\rho), u\rightarrow 0, \rho+\lambda>0.
\end{eqnarray}
The generalized Mellin transform of $h$ is defined by
\begin{eqnarray}
\label{eq:Int-327.2}
M[h; z_1]=\lim_{\varepsilon \rightarrow \, 0+}\int^{\infty}_{0}
u^{z_1-1}h(u) e^{-\varepsilon u} du.
\end{eqnarray}
Then by (\ref{eq:Int-319.2}),
\begin{eqnarray}
\label{eq:Int-328.2}
I_1(c) = \sum^{n-1}_{s=0}c_s M[h;
s+\lambda]c^{-s-\lambda}+\delta_n^1(c),
\end{eqnarray}
where
\begin{eqnarray}
\label{eq:Int-329.2}
\delta_n^1(c)=\lim_{\varepsilon \rightarrow
0+}\int^{\infty}_{0} \overline{\hat{\psi}_n}(u)h(cu)e^{-\varepsilon u} du,
\end{eqnarray}
and, from (\ref{eq:Int-327.2})
\begin{eqnarray*}
M[h(-u); z_1]=\lim_{\varepsilon \rightarrow \, 0+}\int^{\infty}_{0}
u^{z_1-1}h(-u) e^{-\varepsilon u} du.
\end{eqnarray*}
 Hence
\begin{eqnarray}
\label{eq:Int-330.2}
I_2(c) = \sum^{n-1}_{s=0}c_s \left(-1\right)^{s+\lambda+1}M[h(-u);
s+\lambda]c^{-s-\lambda}+\delta_n^2(c),
\end{eqnarray}
where
\begin{eqnarray}
\label{eq:Int-331.2}
\delta_n^2(c)=\lim_{\varepsilon \rightarrow
0+}\int^{\infty}_{0} \overline{\hat{\psi}_n}(-u)h(-cu)e^{-\varepsilon u} du.
\end{eqnarray}
Finally, from (\ref{eq:Int-322.2}), (\ref{eq:Int-328.2}) and (\ref{eq:Int-330.2}) we get the asymptotic expansion:
\begin{eqnarray*}
\left(W_\psi f \right)(b,a) &=& \frac{\sqrt{c}}{2\pi}\biggr\{\sum_{s=0}^{n-1}c_s \left(M[h(u);s+\lambda]+ \left(-1\right)^{s+\lambda+1}M[h(-u);s+\lambda]\right)\\ && \times \; c^{-s-\lambda}+ \delta_n^1(c)+ \delta_n^2(c)\biggr\}
\end{eqnarray*}
Finally, setting $ c=\frac{1}{a} $ we get the asymptotic expansion for small values of a:
\begin{eqnarray}
\label{eq:Int-332.2}
\left(W_\psi f \right)(b,a) \nonumber &=& \frac{1}{2\pi}\biggr\{\sum_{s=0}^{n-1}c_s \left(M[h(u):s+\lambda]+ \left(-1\right)^{s+\lambda+1}M[h(-u):s+\lambda]\right)\\ && \times \; a^{s+\lambda-1/2}+ \delta_n(a)\biggr\},
\end{eqnarray}
where
\begin{eqnarray}
\label{eq:Int-333.2}
\delta_n(a)\nonumber &=&\frac{1}{\sqrt{a}}\lim_{\varepsilon \rightarrow
    0+}\biggr(\int^{\infty}_{0}\overline{\hat{\psi}_n}(u)h(u/a)e^{-\varepsilon u} du \\ && + \;\int^{\infty}_{0}\overline{\hat{\psi}_n}(-u)h(-u/a)e^{-\varepsilon u} du \biggr).
\end{eqnarray}
Using Theorem\ref{thm:ch1sec2-2.1}. we get the following existance theorem for formula (\ref{eq:Int-332.2}).
\begin{thm}
\label{thm:ch1sec2-3.1}
Assume that (i) $ \overline{\hat{\psi}}^{(m)}(u) $ is continuous on $ (-\infty, \infty),$
where $m$ is a nonegative integer; (ii) $ \overline{\hat{\psi}}(u) $ has asymptotic expansion of the form(\ref{eq:Int-325.2}) and the expansion is $m$ times differential; (iii) $h(u)$ satisfies (\ref{eq:Int-324.2}) and (\ref{eq:Int-326.2}) and (iv) as $u \rightarrow \infty \,\,\, u^{-\beta}\overline{\hat{\psi}}^{(m)}(u)= O(u^{-1-\epsilon})$ for $ j= 0,1,2,....m $ and for some $ \epsilon>0 $.
Under these conditions, the result(\ref{eq:Int-332.2}) holds with
\begin{eqnarray}
\label{eq:Int-334.2}
\delta_n(a)\nonumber &=& (-1)^m a^{m-1/2}\int^{\infty}_{-\infty}\overline{\hat{\psi}_n}^{(m)}(u)h^{(-m)}(u/a) e^{-\varepsilon u} du.
\end{eqnarray}
where $n$ is the the smallest positive integer such that $ \lambda + n > m $.
\end{thm}
In the following sections we shall obtain asymptotic expansions for certain special cases of the general wavelet transform.
\section{MORLET WAVELET TRANSFORM}
In this section we shall exploit the following result ~\cite[eq.(12) p.57]{Rainville} for series manipulation
\begin{eqnarray}
\label{eq:Int-335.2}
\sum_{n=0}^\infty\sum_{k=0}^\infty A(k, n)t^{n+2k}= \sum_{n=0}^\infty\sum_{k=0}^{[n/2]}A(k, n-2k)t^{n-2k}.
\end{eqnarray}
We choose
\begin{eqnarray}
\label{eq:Int-336.2}
\psi(t)=\sqrt{2\pi} e^{iu_0 t-t^2/2}.
\end{eqnarray}
Then from ~\cite[p.373]{Debnath}
\begin{eqnarray}
\label{eq:Int-337.2}
\hat{\psi}(u) = \sqrt{2\pi}e^{\frac{-(u-u_0)^2}{2}}.
\end{eqnarray}
Now, using (\ref{eq:Int-335.2}) we can write $\hat{\psi}(u)$ in form of (\ref{eq:Int-325.2})
\begin{eqnarray}
\label{eq:Int-338.2}
\hat\psi(u)
       \nonumber &=& \sqrt{2\pi}e^{-u_0^2/2}e^{u_0 u}e^{-u^2/2} \\
      \nonumber &=& \sqrt{2\pi}e^{-u_0^2/2}\sum_{s=0}^\infty\frac{(u_0u)^s}{s!}\sum_{p=0}^\infty\frac{(-1)^p u^{2p}}{2^p p!}
      \\ \nonumber & =& \sqrt{2\pi}e^{-u_0^2/2}\sum_{s=0}^\infty\sum_{p=0}^\infty\frac{(-1)^pu_0^s u^{s+2p}}{s!2^p p!}
      \\ \nonumber &=& \sqrt{2\pi}e^{-u_0^2/2}\sum_{s=0}^\infty\sum_{p=0}^{[s/2]}\frac{(-1)^p u_0^{s-2p}u^s}{p!(s-2p)! 2^p}
      \\  &=&
      \sum_{s=0}^\infty c_su^s+\hat{\psi}(u),
\end{eqnarray}
where
\begin{eqnarray}
\label{eq:Int-339.2}
c_s=\sqrt{2\pi}e^{-u_0^2/2}\sum_{p=0}^{[s/2]}\frac{(-1)^p u_0^{s-2p}}{p!(s-2p)! 2^p}  .
      \end{eqnarray}
Thus $\hat{\psi}(u)$ possesses asymptotic expansion (\ref{eq:Int-325.2}) with $ \lambda=1 $   and $ c_s $ given by (\ref{eq:Int-337.2}).
Hence, using (\ref{eq:Int-332.2}) and (\ref{eq:Int-337.2}) we get the following asymptotic expansion of $ \left(W_\psi f \right)(b,a) $ for small values of $a$.
\begin{eqnarray}
\label{eq:Int-340.2}
\left(W_\psi f \right)(b,a) \nonumber &=& \frac{1}{2\pi}\biggr\{\sum_{s=0}^{n-1}c_s \left(M[h(u):s+1]+ \left(-1\right)^{s}M[h(-u):s+1]\right)\\ && \times \; a^{s+1/2}+ \delta_n(a)\biggr\},
\end{eqnarray}
where
\begin{eqnarray}
\label{eq:Int-341.2}
\delta_n(a)\nonumber &=&\frac{1}{\sqrt{a}}\lim_{\varepsilon \rightarrow
    0+}\biggr(\int^{\infty}_{0} \hat{\psi}_n(u)h(u/a)e^{-\varepsilon u} du \\ && + \;\int^{\infty}_{0}\hat{\psi}_n(-u)h(-u/a)e^{-\varepsilon u} du \biggr).
\end{eqnarray}
Using Theorem \ref{thm:ch1sec2-3.1}. we get the following existance theorem for formula (\ref{eq:Int-340.2}).
\begin{thm}
\label{thm:ch1sec2-4.1}
Assume that $h(u)$ satisfies (\ref{eq:Int-324.2}) and (\ref{eq:Int-326.2}). Under these conditions, the result (\ref{eq:Int-340.2}) holds with
\begin{eqnarray}
\label{eq:Int-342.2}
\delta_n(a)\nonumber &=& (-1)^m a^{m-1/2}\int^{\infty}_{-\infty}\hat{\psi}_n^{(m)}(u)h^{(-m)}(u/a)e^{-\varepsilon u} du.
\end{eqnarray}
where $n$ is the the smallest positive integer such that $ 1 + n > m $.
\end{thm}
\section{MEXICAN HAT WAVELET TRANSFORM}
In this section we choose
\begin{eqnarray*}
\psi(t)= (1-t^2) e^{-t^2/2}.
\end{eqnarray*}
Then from ~\cite[p. 372]{Debnath}
\begin{eqnarray*}
\hat{\psi}(u)= \sqrt{2\pi} u^2 e^{-u^2/2}.
\end{eqnarray*}
Now,
\begin{eqnarray}
\hat{\psi}(u)\nonumber&=& \sqrt{2\pi} u^2 \sum_{r=0}^\infty \frac{(-1)^r u^{2r}}{2^r r!} \\&=& \nonumber
 \sqrt{2\pi}  \sum_{l=1}^\infty \frac{(-1)^{l-1}u^{2l}}{2^{l-1}(l-1)!}  \\ &=& \sum_{s=0}^\infty c_s u^s,
\end{eqnarray}
where
\begin{eqnarray}
\label{eq:Int-343.2}
c_s = \left\{
        \begin{array}{ll}
           \sqrt{2\pi} \frac{(-1)^{l-1}}{2^{l-1} (l-1)!} & if s=2l, l=1,2,3... \\
             \,\,\,\, 0 & \textrm{otherwise}.
         \end{array}
         \right.
\end{eqnarray}
Then, from (\ref{eq:Int-332.2}) with $ \lambda=1 $ and $ c_s $ given by (\ref{eq:Int-343.2}) yields the asymptotic espansion:
\begin{eqnarray}
\label{eq:Int-344.2}
\left(W_\psi f \right)(b,a) \nonumber &=& \frac{1}{2\pi}\biggr\{\sum_{s=0}^{n-1}c_s \left(M[h(u):s+1]+ \left(-1\right)^{s}M[h(-u):s+1]\right)\\ && \times \; a^{s+1/2}+ \delta_n(a)\biggr\},
\end{eqnarray}
where
\begin{eqnarray}
\label{eq:Int-345.2}
\delta_n(a)\nonumber &=&\frac{1}{\sqrt{a}}\lim_{\varepsilon \rightarrow
    0+}\biggr(\int^{\infty}_{0}\hat{\psi}_n(u)h(u/a)e^{-\varepsilon u} du \\ && + \;\int^{\infty}_{0}\hat{\psi}_n(-u)h(-u/a)e^{-\varepsilon u} du \biggr).
\end{eqnarray}
Using Theorem \ref{thm:ch1sec2-3.1}. we get the following existance theorem for formula (\ref{eq:Int-344.2}).
\begin{thm}
\label{thm:ch1sec2-5.1}
Assume that $h(u)$ satisfies (2.3) and (2.5). Under these conditions, the result (\ref{eq:Int-343.2}) holds with
\begin{eqnarray}
\label{eq:Int-346.2}
\delta_n(a) &=& (-1)^m a^{m-1/2}\int^{\infty}_{-\infty} \hat{\psi}_n^{(m)}(u) h^{(-m)}(u/a)e^{-\varepsilon u} du,
\end{eqnarray}
where $n$ is the the smallest positive integer such that $ 1 + n > m $.
\end{thm}
\section{HAAR WAVELET TRANSFORM}
let us choose\\
\begin{center}
$ \psi(t)=$
$\begin{cases}
\,\,\, 1,\,\,\ 0 \leq t <1/2 & \\ -1 ,\,\, 1/2 \leq t <1 & \\ \,\,\, 0,\,\,\ otherwise,
\end{cases}$\\
\end{center}
Then from ~\cite[p. 368]{Debnath},
\begin{eqnarray}
\label{eq:Int-347.2}
\overline{\hat\psi}(a u)
       \nonumber &=& 4i e^{-iu/2}\frac{sin^2 u/4}{u}\\ \nonumber &=&\frac{i}{u}(1-2 e^{iu/2}+ e^{iu}) \\
      \nonumber &=& \frac{i}{u}\left(1 - 2\sum_{r=0}^\infty \frac{(iu)^r}{2^r r!}+\sum_{r=0}^\infty \frac{(iu)^r}{ r!}\right)
      \\ \nonumber & =& \sum_{r=1}^\infty \frac{i^{r+1}u^{r-1}}{ r!}\left(1-\frac{1}{2^{r-1}}\right)
      \\ \nonumber &=& \sum_{s=0}^\infty \frac{i^{s+2}u^s}{ (s+1)!}\left(1-\frac{1}{2^s}\right)
      \\  &=&
      \sum_{s=0}^\infty c_su^s,
\end{eqnarray}
where
\begin{eqnarray}
\label{eq:Int-348.2}
c_s=\frac{i^{s+2}}{ (s+1)!}\left(1-\frac{1}{2^s}\right).
      \end{eqnarray}
      Then, from (\ref{eq:Int-332.2}) with $ \lambda=1 $ and $ c_s $ given by (\ref{eq:Int-348.2}) we get
\begin{eqnarray}
\label{eq:Int-349.2}
\left(W_\psi f \right)(b,a) \nonumber &=& \frac{1}{2\pi}\biggr\{\sum_{s=0}^{n-1}c_s \left(M[h(u);s+1]+ \left(-1\right)^{s}M[h(-u);s+1]\right)\\ && \times \; a^{s+1/2}+ \delta_n(a)\biggr\},
\end{eqnarray}
where
\begin{eqnarray}
\label{eq:Int-350.2}
\delta_n(a)\nonumber &=&\frac{1}{\sqrt{a}}\lim_{\varepsilon \rightarrow
    0+}\biggr(\int^{\infty}_{0}\overline{\hat{\psi}_n}(u)h(u/a)e^{-\varepsilon u} du \\ && + \;\int^{\infty}_{0}\overline{\hat{\psi}_n}(-u)h(-u/a)e^{-\varepsilon u} du \biggr).
\end{eqnarray}
\section{ASYMPTOTIC EXPANSION FOR SMALL $a$ CONTINUED}
In this section we obtain asymptotic expansion of the wavelet
transform given in the form (\ref{eq:Int-310.2}) when $a \rightarrow 0+$.
Naturally, in this case we have to impose conditions on $f$ and
$\psi$ instead of $\hat{f}$ and $\hat{\psi}.$\\
Now, let us write (\ref{eq:Int-310.2}) in the form:
\begin{eqnarray}
\label{eq:Int-351.2}
(W_\psi f)(b,a)=c^{1/2}\int^{\infty}_{-\infty}f(t+b)
\overline{\psi(ct)}dt,
\end{eqnarray}
where $c=1/a \rightarrow +\infty$ and $b$ is assumed to be a fixed
real number.
\vspace{.5pc} Then setting $g(t) = f(t+b)$ and $h(t) =
\overline{\psi(t)},$ we have
\begin{eqnarray}
\label{eq:Int-352.2}
    (W_\psi f)(b,a) \nonumber &=&c^{1/2}\left[\int^{\infty}_{0}g(t) h(ct) dt +
    \int^{0}_{-\infty}g(t) h(ct) dt\right]\\
    &=&c^{1/2}[I_1+I_2]\;\;\;(say).
    \end{eqnarray}
Assume that $g(t)$ satisfies (\ref{eq:Int-315.2}) and $h(t)$ satisfies (\ref{eq:Int-316.2})
and (\ref{eq:Int-317.2}). Then from (\ref{eq:Int-319.2}) it follows that
\begin{eqnarray}
\label{eq:Int-353.2}
I_1=\sum^{n-1}_{s=0}c_s M[\bar{\psi};
s+\lambda]c^{-s-\lambda}+\delta^{1}_{n}(a),
\end{eqnarray}
where
\begin{eqnarray}
\label{eq:Int-354.2}
\delta^{1}_{n}(a)=\lim_{\varepsilon\rightarrow 0+}\int^{\infty}_{0}
g_n(t) \overline{\psi(t/a)} e^{-\varepsilon t^p}dt;
\end{eqnarray}
and
\begin{eqnarray}
\label{eq:Int-355.2}
I_2=\sum^{n-1}_{s=0}c_s (-1)^{s+\lambda-1}M[\overline{\psi(-t)};
s+\lambda]+\delta^{2}_{n}(a),
\end{eqnarray}
where
\begin{eqnarray}
\label{eq:Int-356.2}
\delta^{2}_{n}(a)=\lim_{\varepsilon\rightarrow 0+}\int^{\infty}_{0}
g_n(t) \overline{\psi(-t/a)} e^{-\varepsilon t^p}dt.
\end{eqnarray}
From (\ref{eq:Int-352.2}), (\ref{eq:Int-353.2}) and (\ref{eq:Int-355.2}) we get
\begin{eqnarray}
\label{eq:Int-357.2}
(W_\psi f)(b, a) \nonumber &=&\sum^{n-1}_{s=0}c_s M[\bar{\psi};
s+\lambda]a^{s+\lambda-1/2}\\ \nonumber
&& + \; \sum^{n-1}_{s=0}c_s (-1)^{s+\lambda-1}M[\overline{\psi(-t)};
s+\lambda]\\ && \times \; a^{s+\lambda-1/2}+\delta_n(a),
\end{eqnarray}
where
\begin{eqnarray}
\label{eq:Int-358.2}
\delta_n(a)\nonumber &=& a^{-1/2}\lim_{\varepsilon \rightarrow 0+}\biggr\{
\int^{\infty}_{0}g_n(t) \overline{\psi(t/a)} e^{-\varepsilon t^p}dt \\  &&
+ \; \int^{\infty}_{0}g_n(-t)
\overline{\psi(- t/a)} e^{-\varepsilon t^p}dt\biggr\}.
\end{eqnarray}
\section{Example}
\label{eq:Int-362.2}
 Let us find again asymptotic
expansion of Morlet wavelet transform for small a , using the above
technique. Here
\[
\psi(t)=e^{i\omega_0 t-t^2/2}.
\]
\vspace{.5pc} Suppose that $g(t) = f(t+b)$ satisfies (\ref{eq:Int-315.2}). Then
from (\ref{eq:Int-353.2}), using formula ~\cite[eq.(21), p.16]{Erdelyi}, we get
\begin{align}
\label{eq:Int-363.2}
    I_1\nonumber &=\sum^{n-1}_{s=0}c_s M[\bar{\psi};
    s+\lambda]c^{-s-\lambda}+\delta^{1}_{n}(a)\\
    \nonumber &=\sum^{n-1}_{s=0}c_s \int^{\infty}_{0}e^{i\omega_0
    t-t^2/2}t^{s+\lambda-1}dt a^{s+\lambda}+\delta^{1}_{n}(a)\\
   &=\sum^{n-1}_{s=0}c_s\Gamma(s+\lambda)e^{-\omega^{2}_{0}/4}D_{-s-\lambda}(-i\omega_0)a^{s+\lambda}+\delta^{1}_{n}(a).
    \end{align}
where $D_\nu(z)$ denoted parabolic cylinder function. Similarly,
\begin{align}
\label{eq:Int-359.2}
I_2
=\sum^{n-1}_{s=0}c_s(-1)^{s+\lambda-1}\Gamma(s+\lambda)e^{-\omega^{2}_{0}/4}D_{-s-\lambda}(i\omega_0)a^{s+\lambda}+\delta^{2}_{n}(a).
   \end{align}
Therefore,
\begin{align}
\label{eq:Int-360.2}
\left(W_\psi \right)(b,
a)\nonumber  =&\sum^{n-1}_{s=0}c_s\Gamma(s+\lambda)e^{-\omega^{2}_{0}/4}(D_{-s-\lambda}(-i\omega_0)\\ &  +(-1)^{s+\lambda-1}D_{-s-\lambda}(i\omega_0))\times a^{s+\lambda-1/2}
+\delta_n(a),
\end{align}
where
\begin{eqnarray}
\label{eq:Int-361.2}
\delta_n(a) \nonumber &=&a^{-1/2}\biggr\{\int^{\infty}_{0}g_n(t)e^{-i\omega_0(t/a)-(t/a)^2/2}dt \\ && + \; \int^{\infty}_{0}g_n(-t)e^{i\omega_0(t/a)-(t/a)^2/2}dt\biggr\}.
\end{eqnarray}

\begin{flushleft}
R S Pathak \\
DST Center for Interdisciplinary Mathematical Sciences \\
Banaras Hindu University \\
Varanasi-221005,India \\
e-mail: ranshankarpathak@yahoo.com \\
\vspace*{0.5cm}
Ashish Pathak \\
Department of Mathemtics and Statistics
\\  Dr. Harisingh Gour Central University
\\   Sagar-470003, India.
\\  e-mail: ashishpathak@dhsgsu.ac.in \\  \
   \end{flushleft}

\thebibliography{00}
\bibitem{Debnath} Debnath, Lokenath; Wavelet Transforms and their Applications,Birkh\"auser (2002).
\bibitem{Erdelyi}Erde'lyi, A.,(Editor) Higher Transcendental functions, vol II, McGram-Hill Book Co., New York (1953).
\bibitem{Rpathak}Pathak, R.S. and Ashish Pathak, Asymptotic expansion of the waveet transform with error term, arXiv:submit/0949812 [math.FA] 4 Apr 2014.
\bibitem{Rainville} Rainvlle, E., Special Functions, macmillan Co. New York (1960).
\bibitem{Wong} Wong, R., Explicit error terms for asymptotic expansion of Mellin convolutions , J. Math. Anal. Appl. 72(1979), 740-756.
\bibitem{Rwong} Wong, R., Asymptotic Approximations of Integrals, Acad. Press, New York (1989).
\bibitem{Aapathak} R S Pathak and Ashish Pathak, Asymptotic Expansions of the Wavelet Transform for Large and Small Values of b, International Journal of Mathematics and Mathematical Sciences, Vol. 2009, 13 page, doi:10.1155/2009/270492.
\\
\end{document}